\documentstyle[12pt]{article}

\makeatletter
   \def\@begintheorem#1#2{\sl \trivlist \item[\hskip \labelsep{\bf #2\ #1}]}
   \def\@opargbegintheorem#1#2#3{\sl \trivlist
            \item[\hskip \labelsep{\bf #1\ #2\ (#3)}]}

   \def\section{\@startsection {section}{1}{\z@}{-3.5ex plus -1ex minus
    -.2ex}{2.3ex plus .2ex}{\large\bf}} %section labels smaller
\makeatother



\makeatletter
\newenvironment{subeqn}{\refstepcounter{subsubsection}
$$}{\leqno{\rm(\thesubsubsection)}$$\global\@ignoretrue}
\makeatother



\newenvironment{prf}[1]{\trivlist
\item[\hskip \labelsep{\it
#1.\hspace*{.3em}}]}{~\hspace{\fill}~$\Box$\endtrivlist}

\newenvironment{proof}{\begin{prf}{\bf Proof}}{\end{prf}}
\let\tempcirc=\circ
\def\circ{\mathord{\raise0.25ex\hbox{$\scriptscriptstyle\tempcirc$}}}

\newcommand{\Bbb}{\bf}

\newcommand{\ZZ}{{\Bbb Z}}
\newcommand{\FF}{{\Bbb F}}
\newcommand{\GL}{{\rm GL}}

\newcommand{\QQ}{{\Bbb Q}}
\newcommand{\Qbar}{{\overline{\QQ}}}

\newcommand{\ld}{\langle}
\newcommand{\rd}{\rangle}
\newcommand{\Aut}{{\rm Aut}}
\newcommand{\Hom}{{\rm Hom}}

\newcommand{\Res}{{\rm Res}}
\newcommand{\nonsplit}{{\mbox{\rm\scriptsize non-split}}}
\newcommand{\jac}{{\rm jac}}
\newcommand{\lto}{\longrightarrow}
\newcommand{\pr}{{\rm pr}}
\newcommand{\ol}{\overline}
\newcommand{\End}{{\rm End}}
\newcommand{\AV}{{\rm AV}}
\newcommand{\QAV}{{\QQ\otimes\AV}}
\newcommand{\calC}{{\cal C}}
\newcommand{\im}{{\rm im}}

\newtheorem{theorem}[subsection]{Theorem.}
\newtheorem{proposition}[subsection]{Proposition.}
\newtheorem{lemma}[subsection]{Lemma.}

\newtheorem{tabel}[subsection]{Table.}

\begin{document}
\title{On a result of Imin Chen.}
\author{Bas Edixhoven}
\maketitle

\section{Introduction, notation and results.}\label{section1}
The aim of this text is to give another proof of a recent result of
Imin Chen, concerning certain identities among zeta functions
of modular curves, or, equivalently, isogenies between products of
jacobians of these curves. I want to thank Imin Chen for pointing
out a mistake in an earlier version of this text.

For $n\geq1$ an integer, let $X(n)_\QQ$ be the modular
curve which is the compactified moduli space (coarse if $n<3$)
of pairs $(E/S/\QQ,\phi)$, where $S$ is a $\QQ$-scheme, $E/S$ is
an elliptic curve and $\phi\colon (\ZZ/n\ZZ)_S^2\to E[n]$ an
isomorphism of group schemes over~$S$. By construction, the group
$\GL_2(\ZZ/n\ZZ)$ acts from the right on $X(n)_\QQ$: an
element $g$ sends $(E/S/\QQ,\phi)$ to $(E/S/\QQ,\phi\circ g)$.
This action induces a left action of the jacobian $J(n)_\QQ$
of~$X(n)_\QQ$.

Let $p$ be a prime number. Let $X$ denote $X(p)_\QQ$ and $G$ the
group $\GL_2(\FF_p)$.
We will consider the following subgroups of $G$: the standard
``maximal torus'' $T$ consisting of diagonal matrices, a non-split
maximal torus $T'$ obtained by choosing an $\FF_p$-basis of
a field $\FF_{p^2}$ of $p^2$ elements, the normalizers $N$ of $T$
and $N'$ of~$T'$. Note that $N/T$ and $N'/T'$ are both of order~2.
Finally, let $B_+$ and $B_-$ denote the two Borel subgroups
containing $T$; $B_+$ is the subgroup of upper triangular matrices
and $B_-$ the one of lower triangular matrices.

The quotients of $X$ by some of these subgroups have the following
interpretations. The quotient $X/T'$ is usually denoted
$X(p)_\nonsplit$. The constructions
\begin{subeqn}\label{eqn1.0.1}
(E/S,\phi) \mapsto (E/S,\ld \phi(1,0)\rd), \qquad
(E/S,\phi) \mapsto (E/S,\ld \phi(0,1)\rd)
\end{subeqn}
induce isomorphisms
\begin{subeqn}\label{eqn1.0.2}
X/B_+\;\;\tilde{\lto}\;\; X_0(p)_\QQ, \qquad
X/B_-\;\;\tilde{\lto}\;\; X_0(p)_\QQ
\end{subeqn}
The construction
\begin{subeqn}\label{eqn1.0.3}
(E/S,\phi) \mapsto (E_1/S,\ker(\phi_2\circ\phi_1^*))
\end{subeqn}
where $\phi_1\colon E\to E_1$ (resp. $\phi_2\colon E\to E_2$) is the
isogeny whose kernel is the subgroup scheme generated by $\phi(1,0)$
(resp. $\phi(0,1)$), induces an isomorphism
\begin{subeqn}\label{eqn1.0.4}
X/T\;\;\tilde{\lto}\;\;X_0(p^2)_\QQ
\end{subeqn}
Under this isomorphism the Atkin-Lehner
involution $w_{p^2}$ of $X_0(p^2)_\QQ$ corresponds to the non-trivial
element of~$N/T$; the two maps $X/T\to X/B_+$ and $X/T\to X/B_-$
correspond to the two standard degeneracy maps from $X_0(p^2)_\QQ$
to~$X_0(p)_\QQ$.

The result of Chen is the following, see \cite[Theorem~1 and \S10]{Chen}.

\begin{theorem}\label{theorem1.1} {\bf (Chen)}
The jacobian of $X_0(p^2)_\QQ$ is isogeneous to the product of the
jacobian of $X(p)_\nonsplit$ by the square of the jacobian of
$X_0(p)_\QQ$. The jacobian of $X_0(p^2)_\QQ/\ld w_{p^2}\rd$
is isogeneous to the product of the jacobian of $X/N'$ by the
jacobian of~$X_0(p)_\QQ$.
\end{theorem}
The proof given by Chen is to show that the traces of the Hecke
operators $T_n$ ($n$ prime to $p$) on the jacobians in the theorem
satisfy the identities required to conclude by the Eichler--Shimura
relations and Faltings's isogeny theorem that one has the desired
isogenies. We will prove a generalization of Theorem~\ref{theorem1.1}
using only the representation theory of~$G$ and some elementary properties
of abelian varieties.

For a field $k$, let $\AV(k)$ denote the category of abelian varieties
over~$k$. Let $\QAV(k)$ denote the category of abelian varieties over $k$
``up to isogeny'', i.e., its objects are those of $\AV(k)$ and for two
objects $A$ and $B$ one has
$\Hom_{\QAV(k)}(A,B) = \QQ\otimes\Hom_{\AV(k)}(A,B)$. For $A$ an abelian
variety over $k$ we denote by $\QQ\otimes A$ the corresponding object
of~$\QAV(k)$. By construction, $A$ and $B$ are isogeneous if and only
if $\QQ\otimes A$ and $\QQ\otimes B$ are isomorphic. The categories
$\QAV(k)$ are $\QQ$-linear, semi-simple and abelian.

Recall (e.g., see \cite[\S1]{Scholl1}), that an additive category $\calC$
is called pseudoabelian if for every object $M$ of $\calC$ every
idempotent $f$ in $\End(M)$ has a kernel (or, equivalently, an image).
If $\calC$ is additive, pseudoabelian and $f$ in $\End(M)$ is an
idempotent in $\calC$, then the natural morphism from $\im(f)\oplus\ker(f)$
to $M$ is an isomorphism. The categories $\QAV(k)$ are clearly additive
and pseudoabelian.

For each subgroup $H$ of $G$ we define
\begin{subeqn}\label{eqn1.1.1}
\pr_H := \frac{1}{|H|}\sum_{h\in H} h \in \QQ[G]
\end{subeqn}
Hence $\pr_H$ is the idempotent of $\QQ[G]$ that projects on the
$H$-invariants. For two subgroups $H_1$ and $H_2$ of $G$ such that
$\ld H_1\cup H_2\rd = H_1H_2$, one has
$\pr_{H_1}\pr_{H_2}=\pr_{\ld H_1\cup H_2\rd}$. For $H$ a subgroup
and $g$ in $G$ one has $g\pr_Hg^{-1}=\pr_{gHg^{-1}}$, hence $\pr_H$ is
a central idempotent if and only if $H$ is a normal subgroup.

For each irreducible representation $V$ of $G$ over $\QQ$ let $e_V$
be the corresponding central idempotent in~$\QQ[G]$ which projects on
the $V$-isotypical part. If $V$ is absolutely irreducible, of dimension~$d$
and with character $\chi$, one has:
\begin{subeqn}\label{eqn1.1.2}
e_V := \frac{d}{|G|}\sum_{g\in G} \chi(g^{-1})g
\end{subeqn}
We will use only one idempotent of the form $e_V$, namely, with $V$ the
representation with character $\pi^-(1)$ (see Table~\ref{table2.1}). This
representation is the $p$-dimensional irreducible subrepresentation of
the induction of the trivial representation from $B_+$ to~$G$. It is
clearly absolutely irreducible and it exists over~$\QQ$.

Let us for the moment admit the following proposition, whose proof
will be given in the next section.
\begin{proposition}\label{proposition1.2}
Suppose that $p\neq2$. The elements
$\pr_{T'}(1-\pr_G)$ and $\pr_T(1-e_{\pi^-(1)})(1-\pr_G)$
of the ring $\QQ[G]$ are conjugate idempotents. Likewise, the elements
$(\pr_{N'}+\pr_{B_+})(1-\pr_G)$ and $\pr_{N}(1-\pr_G)$ are conjugate
idempotents.
\end{proposition}
Our generalization of Chen's result is simply the following direct
consequence of Proposition~\ref{proposition1.2}.
\begin{theorem}\label{theorem1.3}
Suppose that $p\neq2$. Take elements $u$ and $v$ of $\QQ[G]^*$ such that
\begin{eqnarray*}
u\pr_{T'}(1-\pr_G)u^{-1} & = & \pr_T(1-e_{\pi^-(1)})(1-\pr_G)\\
v(\pr_{N'}+\pr_{B_+})(1-\pr_G)v^{-1} & = & \pr_{N}(1-\pr_G)
\end{eqnarray*}
Let $\calC$ be a $\QQ$-linear pseudoabelian additive category. Let $M$ be an
object of $\calC$ with an action by the group~$G$; this gives a morphism
of rings $\QQ[G]\to\End(M)$. Then $u$ induces an isomorphism
$$
\pr_{T'}(1-\pr_G)M \;\;\tilde{\lto}\;\; \pr_T(1-e_{\pi^-(1)})(1-\pr_G)M
$$
Likewise, $v$ induces an isomorphism
$$
\pr_{N'}(1-\pr_G)M\oplus\pr_{B_+}(1-\pr_G)M \;\;\tilde{\lto}\;\;
\pr_N(1-\pr_G)M
$$
\end{theorem}
To see that Theorem~\ref{theorem1.1} is a special case, apply
Theorem~\ref{theorem1.3} to $\calC:=\QAV(\QQ)$ and take
$M=\QQ\otimes\jac(X)$, with $\jac(X)$ the jacobian of~$X$. For any subgroup
$H$ of $G$ one then has $\pr_HM=\QQ\otimes\jac(X/H)$. In this case
$\pr_G$ acts as zero on $M$, since $X/G$ has genus zero.
The idempotent $e_{\pi^-(1)}$, acting on
$\QQ\otimes\jac(X/T)=\QQ\otimes J_0(p^2)$, projects on the old part, which is
a product of two copies of $\QQ\otimes J_0(p)$ (one way to see this is to
note that the space of $T$-invariants in the representation corresponding
to $\pi^-(1)$ is the direct sum of the two $1$-dimensional spaces of $B_+$
and $B_-$-invariants).
One also has to
use the interpretations of the $X/H$ as explained in
the beginning of this section. For the case $p=2$, note that $X(2)_\QQ$ has
genus zero.

\section{The proof of Proposition~1.2.}\label{section2}
The notation is as in the previous section, in particular, $G=\GL_2(\FF_p)$.
We suppose that $p\neq2$.
We will need to do some calculations involving the irreducible characters
of~$G$, so for convenience of the reader and to fix the notation, we
include its character table, taken from~\cite{Cartier}:

\begin{tabel}{The character table of $G$.}\label{table2.1}
$$
\renewcommand{\arraystretch}{1.5}
\begin{array}{|l||c|c|c|c|}
\hline
 \mbox{\rm conjugacy class of} & \pi(\alpha,\beta),\;\;\alpha\neq\beta &
 \pi(\Lambda),\;\;\Lambda^p\neq\Lambda & \alpha\circ\det & \pi^-(\alpha) \\
\hline \hline
\bigl({x\atop 0}{0\atop x}\bigr)\;\;x\in\FF_p^* &
(p+1)\alpha(x)\beta(x) & (p-1)\Lambda(x) & \alpha(x)^2 & p\alpha(x)^2 \\
\hline
\bigl({x\atop 0}{0\atop y}\bigr)\;\;x,y\in\FF_p^*\;\;x\neq y &
\alpha(x)\beta(y)+\alpha(y)\beta(x) & 0 & \alpha(x)\alpha(y) &
\alpha(x)\alpha(y) \\
\hline
\bigl({x\atop 0}{1\atop x}\bigr)\;\;x\in\FF_p^* &
\alpha(x)\beta(x) & -\Lambda(x) & \alpha(x)^2 & 0 \\
\hline
\bigl({z\atop 0}{0\atop z^p}\bigr)\;\;z\in\FF_{p^2}^*\;\;z^p\neq z &
0 & -\Lambda(z)-\Lambda(z^p) & \alpha(z^{p+1}) & -\alpha(z^{p+1}) \\
\hline
\end{array}
$$
\end{tabel}
In this table $\alpha$ and $\beta$ denote characters $\FF_p^*\to\Qbar^*$ and
$\Lambda$ denotes a character $\FF_{p^2}^*\to\Qbar^*$. For each effective
character $\chi$ of $G$ we denote by $V_\chi$ some $\Qbar[G]$-module with
character~$\chi$. For each irreducible $\chi$ and each of the subgroups
$H\subset G$ mentioned at the beginning of \S\ref{section1}, we will need to
know the dimension $\dim(V_\chi^H)$ of the set of $H$-invariants in~$V_\chi$.
These dimensions are given in the following table, in which
$\delta(x,y)$ denotes the Kronecker symbol, i.e., $\delta(x,y)=1$ if $x=y$
and $\delta(x,y)=0$ otherwise.
\begin{tabel}{The dimensions of the spaces $V_\chi^H$. }\label{tabel2.2}
$$
\renewcommand{\arraystretch}{1.5}
\begin{array}{|l||c|c|c|c|}
\hline
 & \pi(\alpha,\beta) & \pi^-(\alpha) & \alpha\circ\det & \pi(\Lambda) \\
\hline \hline
T & \delta(\alpha\beta,1) & \delta(\alpha,1)+\delta(\alpha^2,1) &
\delta(\alpha,1) & \delta(\Lambda^{p+1},1) \\
\hline
N & \delta(\alpha(-1),1)\delta(\alpha\beta,1) &
\delta(\alpha(-1),1)\delta(\alpha^2,1) & \delta(\alpha,1) &
\delta(\Lambda^{p+1},1)-\delta(\Lambda^{(p+1)/2},1) \\
\hline
T' & \delta(\alpha\beta,1) & -\delta(\alpha,1)+\delta(\alpha^2,1) &
\delta(\alpha,1) & \delta(\Lambda^{p+1},1) \\
\hline
N' & \delta(\alpha(-1),1)\delta(\alpha\beta,1) &
-\delta(\alpha,1)+\delta(\alpha(-1),1)\delta(\alpha^2,1) &
\delta(\alpha,1) &
\delta(\Lambda^{p+1},1)-\delta(\Lambda^{(p+1)/2},1) \\
\hline
B & 0 & \delta(\alpha,1) & \delta(\alpha,1) & 0 \\
\hline
\end{array}
$$
\end{tabel}
We will not give the computation of this table in detail, since it is
a straightforward application of the theory of representations of finite
groups, see for example~\cite{Serre1}. As an example, let us do the
case $\chi=\pi(\Lambda)$ and $H=N'$ (the other computations are in fact
easier). The group $N'$ can be identified with the subgroup of
$\GL_{\FF_p}(\FF_{p^2})$ generated by $\FF_{p^2}^*$ and $\sigma$, where
$\sigma$ is the automorphism of order two of~$\FF_{p^2}$. Then $N'$ is
the disjoint union of $T'=\FF_{p^2}^*$ and $T'\sigma$. The conjugacy
class in $G$ of $z\in T'$ is the conjugacy class
of~$({z\atop0}{0\atop z^p})$. The conjugacy class of $z\sigma$ is
the one of~$({z^{(p+1)/2}\atop 0}{0\atop -z^{(p+1)/2}})$.
One has:
\begin{subeqn}\label{eqn2.2.1}
\dim(V_\chi^H) = \dim \Hom_H(\Res^G_H(V_\chi),\Qbar) =
\frac{1}{|H|}\sum_{g\in H}\chi(g)
\end{subeqn}
The sum over the elements of $T'$ can be written as:
\begin{subeqn}\label{eqn2.2.2}
-\sum_z\left(\Lambda(z)+\Lambda(z^p)\right) + (p+1)\sum_x\Lambda(x)
\end{subeqn}
In this sum, $z$ runs through $\FF_{p^2}^*$ and $x$ through $\FF_p^*$.
The first of the two terms of (\ref{eqn2.2.2}) gives zero, the second
contributes ${1\over2}\delta(\Lambda^{p+1},1)$ to~$\dim(V_\chi^H)$.
The sum over the elements of $T'\sigma$ can be written as
\begin{subeqn}\label{eqn2.2.3}
\sum_{z^{(p+1)/2}\in\FF_p}
\left(\Lambda(z^{(p+1)/2})+\Lambda(-z^{(p+1)/2})\right) -
\sum_z\left(\Lambda(z^{(p+1)/2})+\Lambda(-z^{(p+1)/2})\right)
\end{subeqn}
The first of the two terms of (\ref{eqn2.2.3}) contributes
${1\over2}\delta(\Lambda^{p+1},1)$ to~$\dim(V_\chi^H)$ and the second
term contributes~$-\delta(\Lambda^{(p+1)/2},1)$. This completes the
computation of~$\dim(V_\chi^H)$.

As promised, we will now give a proof of Proposition~\ref{proposition1.2}.
In fact, that proposition is a direct consequence of the following one.

\begin{proposition}\label{proposition2.3}
Define $\ol{\QQ[G]}:=\QQ[G]/(\pr_G)$ and denote the projection
$\QQ[G]\to\ol{\QQ[G]}$ by $u\mapsto\ol{u}$. Then the elements
$\ol{\pr_{T'}}$ and $\ol{\pr_T}(1-\ol{e_{\pi^-(1)}})$
of the ring $\ol{\QQ[G]}$ are conjugate idempotents. Likewise, the elements
$\ol{\pr_{N'}}+\ol{\pr_{B_+}}$ and $\ol{\pr_{N}}$ are conjugate idempotents.
\end{proposition}
\begin{proof}
Consider the first statement. Both elements are clearly idempotents.
The $\QQ$-algebra $\QQ[G]$ is a product of matrix algebras over division
rings. Using Table~\ref{tabel2.2} one verifies that the two elements in
question generate, in each factor, two left ideals of the same dimension
over $\QQ$ (actually, one verifies this after extension of scalars to
$\Qbar$). Lemma~\ref{lemma2.4} then implies that the two elements are
conjugates.

The proof of the second statement is almost the same. The element
$\ol{\pr_{N'}}+\ol{\pr_{B_+}}$ is an idempotent because $T'B_+=G$. The rest
of the proof runs as before.
\end{proof}

\begin{lemma}\label{lemma2.4}
Let $\Delta$ be a division ring. Let $0\leq k\leq n$ be integers. Then
the group $\GL_n(\Delta)$ acts transitively (by conjugation) on the set
of idempotents of rank $k$ in ${\rm M}_n(\Delta)$.
\end{lemma}
\begin{proof}
Consider the right $\Delta$-module $\Delta^n$. Then ${\rm M}_n(\Delta)$
can be viewed as $\End_\Delta(\Delta^n)$. The map that associates to an
idempotent of rank $k$ its kernel and image is a bijection between the set
of such idempotents and the set of pairs of $\Delta$-submodules $(V_1,V_2)$
such that $\dim_\Delta(V_2)=k$ and $\Delta^n=V_1\oplus V_2$. One verifies
easily that $\Aut_\Delta(\Delta^n)$ acts transitively on the set of such
pairs.
\end{proof}

\end{document}